\documentclass{amsart}
\usepackage{amsfonts,amscd, amssymb}

\newtheorem{theorem}{Theorem}[section]
\newtheorem{lemma}[theorem]{Lemma}
\newtheorem{corollary}[theorem]{Corollary}
\newtheorem{proposition}[theorem]{Proposition}
\theoremstyle{remark}

\theoremstyle{definition}
\newtheorem{definition}[theorem]{Definition}

\numberwithin{equation}{section} \makeatother

\DeclareMathOperator{\Cdb}{{\mathbb C}}

\DeclareMathOperator{\Ddb}{{\mathbb D}}

\begin{document}

\title[Beurling theorem for
noncommutative $L^p$]{A Beurling theorem for
noncommutative $L^p$}


\author{David P. Blecher}
\address{Department of Mathematics, University of Houston, Houston, TX
77204-3008}
\email[David P. Blecher]{dblecher@math.uh.edu}
 \author{Louis E. Labuschagne}
\address{Department of Mathematical Sciences, P.O. Box
392, 0003 UNISA, South Africa}
\email{labusle@unisa.ac.za}
\thanks{*Blecher is partially supported by a grant from
the National Science Foundation.
Labuschagne is  partially supported by the National Research
Foundation.}
\thanks{The results contained in this
paper in the case $p \leq 2$ were first presented 
at the West Coast Operator Algebras Seminar Seminar, Banff
International Research Station, September 2005.}
\begin{abstract}
We extend Beurling's invariant subspace 
theorem, by characterizing subspaces $K$ of 
the noncommutative $L^p$ spaces which are invariant with
respect to Arveson's maximal subdiagonal algebras,
sometimes known as noncommutative $H^\infty$.  
It is significant that a certain subspace, and a certain quotient,
of $K$ are $L^p({\mathcal D})$-modules in the recent sense
of Junge and Sherman, and therefore have a nice
decomposition into cyclic submodules.       
We also give general inner-outer factorization formulae for elements in the 
noncommutative $L^p$.  
These facts generalize the classical ones, and should
be useful in the future development of noncommutative $H^p$ theory.
In addition, these results characterize maximal subdiagonal algebras.
\end{abstract}

\maketitle

\section{Introduction}

The starting point of this paper is Beurling's invariant
subspace theorem,
stating that a certain class of
invariant subspaces of $L^2$ may be characterized as exactly those spaces
of the form $u H^2$, where $u$ is a unimodular function, and
$H^2$ is the Hilbertian Hardy space.  Many generalizations
of this theorem have appeared over the decades (e.g.\ \cite{Hal}), and our
paper is concerned with generalizations appropriate to
Arveson's noncommutative generalization of the Hardy
spaces of the disk.

Throughout this paper, $M$ is a finite von Neumann algebra
possessing a faithful normal tracial state $\tau$, and
$A$ is a tracial subalgebra of $M$.   That is,
$A$ is a weak* closed unital
subalgebra $A$ of $M$ for which the restriction to $A$ of the unique faithful normal
conditional expectation $\Phi$
from $M$ onto ${\mathcal D} \overset{{\rm def}}{=} A \cap A^*$
satisfying $\tau = \tau \circ \Phi$, is a homomorphism.
Here $A^*$ denotes the set of adjoints of elements
in $A$.  Two simple examples which may help the reader
understand the setting, are 1)  the subalgebra of the $n \times n$ matrices
consisting of the upper
triangular matrices, and 2) the classical $H^\infty$ space of
bounded analytic functions on the disk (here
$\tau = \Phi$ is just the Haar integral on $L^\infty(\mathbb{T})$).
A tracial subalgebra of $M$ is called {\em maximal subdiagonal}
if $A  + A^*$ is weak* dense in $M$.   This form of the definition
is due to Arveson \cite{AIOA} and Exel \cite{E}.   Sometimes such
$A$ is called a {\em noncommutative $H^\infty$}.
If $M$ is commutative and ${\mathcal D}$ is one dimensional
then $A$ is a {\em weak* Dirichlet algebra} \cite{SW}.
In an earlier paper we gave a list of many disparate looking
conditions which a tracial algebra $A$ might satisfy, which turn out to be
equivalent to each other, and equivalent to $A$ being
maximal subdiagonal \cite{BL2}. 

We write $L^p(M)$ for $L^p(M, \tau)$, the
noncommutative $L^p$-space associated to the pair $(M,\tau)$,
in the sense of e.g.\ Nelson \cite{Nel} as
a certain space of operators affiliated to $M$.
 In the present paper we study the structure of
(right) $A$-invariant subspaces of $L^p(M)$, and give general
inner-outer factorization formulae of Beurling-Nevanlinna type for
elements of $L^p(M)$.
These results generalize important classical results
(e.g.\ see references in \cite{SW}), and should
be useful in the future development of the 
noncommutative $H^p$ theory.   They also
constitute a natural occurrence of the `$L^p$-modules' 
and `$L^p$-column sums' due recently to Junge and Sherman \cite{JS}.
In addition, our results characterize maximal subdiagonal algebras,
allowing us to supplement the list given in our earlier paper of
criteria equivalent to maximal subdiagonality.

We will use the notation of \cite{BL2}, to which the reader is
also referred for further explanations and details.
For any set $\mathcal{S}$ of operators, 
$\mathcal{S}^*$ will denote the set of adjoints of
elements in $\mathcal{S}$.
We write $[{\mathcal S}]_p$ for the closure of 
${\mathcal S}$ in $L^p(M)$.  If $A$ is a
 maximal subdiagonal algebra, then
$[A]_p$ is often called a {\em
noncommutative Hardy space}, and written as $H^p$.
We write $A_\infty$ for $[A]_2 \cap
M$ (which equals $A$ if $A$ is
subdiagonal), and $A_0$ for $A \cap {\rm Ker}(\Phi)$.  We say that
$A$ satisfies $L^2$-density, if $A + A^*$ is $L^2$-dense in $L^2(M)$.
Clearly $A$ satisfies $L^2$-density if and only if the same is true
of $A_\infty$.   We say that
$A$ has the {\em unique normal state extension
property} if whenever $g \in
L^1(M)_+$ with $\tau(g A_0) = 0$, then $g \in
L^1({\mathcal D})$.     All maximal subdiagonal algebras 
have these latter two properties.
In fact, it is shown in \cite{BL2} that maximal subdiagonal
algebras are exactly the tracial algebras possessing these two properties.

We recall that a {\em (right) invariant subspace} of $L^p(M)$, is
a closed subspace $K$ of $L^p(M)$ such that $K A \subset K$.
For consistency,
we will not consider left invariant subspaces at all, leaving the reader to verify
that entirely symmetric results pertain in the left invariant case.
An invariant subspace is called {\em simply invariant} if in addition
the closure of $K A_0$ is properly contained in $K$.   
It is the latter class of subspaces to which the generalized Beurling theorem
 for e.g.\ weak* Dirichlet algebras applies.     
In the literature there are several invariant subspace theorems, 
inspired by the Beurling result and its classical extensions, and
associated factorization
results, for maximal subdiagonal algebras (see e.g.\ \cite{Kam,MMS,MW,N,PX,Sai,Zs}).
We mention just two which we shall use: 
Saito showed in \cite{Sai} that any $A$-invariant subspace of
$L^p(M)$ is the closure of the bounded elements which it contains.  Nakazi and Watatani
showed in \cite{N} that in the case that the center
of $M$ contains the center of ${\mathcal D}$, every `type 1' (defined below)
invariant subspace of $L^2(M)$ is of the form $u [A]_2$ for a
partial isometry $u$.  Compelling examples of invariant subspaces exhibiting
interesting structure may be found in \cite{MMS,Zs}.

If $K$ is a right $A$-invariant subspace of $L^2(M)$, we
follow Nakazi and Watatani \cite{N}, defining
the {\em right wandering subspace} of $K$ to
be the space $W = K \ominus [K A_0]_2$;
and we say that $K$ is {\em type 1} if
$W$ generates $K$ as an $A$-module (that is,
$K = [W A]_2$).  
We will say that $K$ is
{\em type 2} if $W = (0)$.  (The last
notation conflicts with that of \cite{N}, where 
this
class of subspaces is decomposed 
into two further subclasses which they call
type II and type III.) 
If $p \neq 2$ one may define the {\em wandering quotient} 
to be $K/[K A_0]_p$, and say that $K$ is
type 2 if this is trivial.   In our paper  
if $p = \infty$, we take $[\cdot]_p$  to be 
the weak* closure.   It turns out that
the wandering quotient is an $L^p({\mathcal D})$-module
in the sense of \cite{JS}, and it is isometric to
a canonically defined subspace of $K$ which can be called
the {\em right wandering subspace of} $K$.
We say that $K$ is type 1 if this subspace generates $K$
as an $A$-module.  
If $1 \leq p < 2$ (resp.\ $p > 2$) then we will show that
$K$ is type 1 iff $K \cap L^2(M)$ (resp.\ $[K]_2$)
is type 1 in the sense of the $L^2$ case above.

In the classical case, or more generally whenever ${\mathcal D}$
is one dimensional, there is a dichotomy:
in this case type 1 is the same thing as being simply invariant,
and any invariant subspace which is not type 1 is type 2.
 In the general case, being type 2 is the same as 
being not simply invariant; and any nontrivial type 1 subspace is 
simply invariant.  However the `simply invariant' 
condition no longer plays a very significant
role for us.  Moreover, there is no longer a dichotomy between
types 1 and 2.  Instead, there is a direct sum decomposition.
We use the {\em column $L^p$-sum} recently studied by Junge and 
Sherman \cite{JS} to investigate this:  If $X$ is a subspace of $L^p(M)$, 
and if  $\{ X_i : i \in I \}$ is
a collection of  subspaces of $X$, which together densely
span $X$, with the
property that $X_i^* X_j = \{ 0 \}$ if $i \neq j$, then
we say that $X$ is the {\em internal column
$L^p$-sum} $\oplus^{col}_i \, X_i$.  If 
$p = \infty$ we also assume that $X$ and $X_i$ are weak* closed,
and the word `densely' above is taken with regard to the
weak* topology.     
Our main result, which builds on earlier ideas from \cite{N,MMS},
is as follows:

\begin{theorem} \label{main}  If $A$ is a maximal subdiagonal
subalgebra of $M$, if $1 \leq p \leq \infty$,
and if $K$ is a closed 
(indeed weak* closed, if $p = \infty$) right  $A$-invariant subspace of $L^p(M)$,
then: \begin{itemize}
\item [(1)]   $K$ may be written uniquely as
an (internal) $L^p$-column sum $K_1 \oplus^{col} K_{2}$
of a type 1 and a type 2 invariant subspace of $L^p(M)$, respectively.
  \item [(2)]  If
$K  \neq (0)$ then $K$ is type 1 if and only if
$K = \oplus_i^{col} \, u_i \, H^p$, for $u_i$  partial isometries
with mutually orthogonal ranges and $|u_i| \in {\mathcal D}$.
\item [(3)]  
The right wandering subspace $W$ of $K$
is a $L^p({\mathcal D})$-module in the sense of Junge and Sherman, and in
particular $W^* W \subset L^{p/2}({\mathcal D})$. 
  \end{itemize}
Conversely, if $A$ is a tracial subalgebra of $M$ such that 
every $A$-invariant subspace of $L^2(M)$
satisfies {\rm (1)} and {\rm (2)} (resp.\ {\rm (1)} and {\rm (3)}), 
then $A$ is maximal subdiagonal.
  \end{theorem}

Note that this 
theorem immediately implies the classical Beurling theorem,
and its generalization to 
simply invariant subspaces of weak* Dirichlet algebras.
Indeed if ${\mathcal D}$ is one-dimensional,
and if $K$ is a simply invariant subspace of $L^p$,
then $K$ is not type 2,
and so
$K_1 \neq (0)$.  Thus there
is a nonzero partial isometry $u$ with $|u| \in {\mathcal D}
= \Cdb 1$.  Thus $u^* u = 1$, and since $M$
is finite $u$ is a unitary in $M$.  Thus $K_1 = u H^p$.
Since $K_1^* K_{2} = (0)$, we have $K_{2} = (0)$.  Thus $K = u H^p$, as desired.

An element $\xi$ of a ${\mathcal D}$-module is called {\em separating}
if the map $d \mapsto \xi d$ is one-to-one on ${\mathcal D}$; and 
{\em cyclic} if $\xi {\mathcal D}$ is dense.  If $p = \infty$ we 
mean dense in the weak* topology.
 
\begin{corollary} \label{stand}  If $A$ and $K$ are as in
Theorem {\rm \ref{main}}, 
then $K$ is of the form $u H^p$ for a unitary $u \in M$,
if and only if the right wandering subspace of $K$ is a `standard' 
representation of ${\mathcal D}$, that is, it has a nonzero
separating and cyclic vector for the right
 action of ${\mathcal D}$.  
This is equivalent to the right wandering quotient
having a separating and cyclic vector.      
\end{corollary}

We say that $f \in L^2(M)$ is
{\em Beurling-Nevanlinna factorizable} (or {\em BN-factorizable}),
if $f = uh$, for a unitary $u$ in $M$ and an $h$ with
$[h A]_2 = [A]_2$.  An $h$ with
$[h A]_2 = [A]_2$ will be called  {\em outer},  as in 
the classical theory.  The content of the assertion in this definition
is that $h \in [A]_2$, and $1 \in [h A]_2$.
The unitary $u$ here is called {\em inner}.
 We will say that $f \in L^2(M)$ is {\em partially  BN-factorizable}
if $f = uh$, for a partial isometry $u$ in $M$ with
$u^* u \in {\mathcal D}$, and an $h$ with
$[h A]_2 = (u^* u) [A]_2$.   
That is,  $h = (u^* u) h \in [A]_2$, and $u^* u \in [h A]_2$.
Our invariant subspace theorem above leads immediately, 
as in the classical case,
 to Beurling-Nevanlinna factorization results,
as we shall see in Section 3.  It will be clear in Section 4
that these results extend to $L^p$ for $p \neq 2$.

It is worth noting that if we assume that our tracial algebra is
{\em antisymmetric}
(that is, ${\mathcal D}$ is one-dimensional, which forces $\Phi(\cdot) = \tau(\cdot)1$), then
almost all of the classical results about generalized Hardy spaces 
found in \cite{Sr,SW}, for example,
and their proofs, seem to transfer almost verbatim
and without difficulty. 
To illustrate this, we end this introduction with 
a special case of our main results.
First we will need a simple lemma:

\begin{lemma}  \label{aeqainf}  Let $A$ be a tracial subalgebra of $M$
which satisfies $[A]_1 = \{ x \in L^1(M) : \tau(xA_0) = 0 \}$.
Then $A = A_\infty$, and $A = \{ f \in M : \tau(f A_0) = 0 \}$.
\end{lemma}

\begin{proof}  This follows just as in \cite[Lemma, p.\ 816]{Sr}.
Following that proof we find an $f \in L^1(M)$ with
$\tau(f A) = 0$.  By our hypothesis, $f \in [A]_1$.  For $d \in {\mathcal D}$ we have
$$0 = \tau(a f) = \tau(\Phi(f d)) = \tau(\Phi(f) d) .$$
It follows that $\Phi(f) = 0$, and so $f \in [A_0]_1$.
As in the last cited reference we conclude that $[A]_1 \cap M = A$,
which implies the last assertion of the lemma.  Also,
$A_\infty = [A]_2  \cap M \subset [A]_1 \cap M = A$, giving
the other assertion.
\end{proof}

\begin{proposition}  \label{rest} Let $A$ be a tracial
subalgebra of $M$.  Consider the conditions:
\begin{itemize} \item [{\rm (a)}]
 Every simply right invariant subspace of $L^2(M)$
is of the
form $u [A]_2$, for a unitary $u$ in $M$,
\item [{\rm (b)}]
Whenever  $f \in L^2(M)$ with $f \notin [f A_0]_2$,
then $f$ is BN-factorizable,
\item [(c)]  $A$ is maximal subdiagonal.
\end{itemize}
Then  {\rm (a)} $\Rightarrow$ {\rm (b)} $\Rightarrow$
{\rm (c)}. If $A$ is antisymmetric, then the conditions
are all equivalent.
\end{proposition}

\begin{proof}
 That {\rm (a)} $\Rightarrow$ {\rm (b)}
follows just as in \cite{SW}, for example, with
only minor modifications.  One needs to use 
the fact that left invertibility implies
invertibility in a finite von Neumann algebra.
Similarly for the implication {\rm (c)} $\Rightarrow$ {\rm (a)} 
in the antisymmetric case
(see also \cite{Kam}), or note that this is proved in the paragraph
after Theorem \ref{main}.
Similarly, supposing {\rm (b)}, it also follows just as in \cite{SW},
that $A_\infty$ is maximal subdiagonal.
We next claim that  $f \in L^1(M)$ and $\tau(f A_0) = 0$,
iff $f \in [A]_1$.  This may be proved as in
\cite[Corollary 2.3]{Sr} (note that \cite[Theorem 2]{Sr} follows
easily from condition {\rm (b)}).
By Lemma \ref{aeqainf},
$A_\infty = A$.  So $A$ is maximal subdiagonal.
\end{proof}

If $A$ is not
antisymmetric, then (c) in this last Proposition need not imply (a) or (b).
This may be seen by considering the example of $M$ a
two-dimensional von Neumann algebra, and $A = M$.  
Also,  for the upper triangular matrices both (a) and (b) fail.
In fact, it
is not hard to see that if ${\mathcal
D}$ is not a factor, then the conditions cannot be equivalent.
Below, we will find the appropriate generalizations of the statements of the
simply invariant subspace theorem, and the Beurling-Nevanlinna
factorization result, and we will show that the new statements are
each equivalent to maximal subdiagonality.  

\section{Invariant subspaces of $L^2(M)$}

We begin with some general observations about the structure of
invariant subspaces.  There is a nontrivial overlap between (5) below, 
and results in \cite{N} (see particularly Lemma 2.4 and Theorem 2.14 there).

 \begin{theorem} \label{inv1}
Let $A$ be a tracial algebra.    \begin{itemize}
\item [(1)]   Suppose that $X$ is a subspace of $L^2(M)$
of the form $X = Z \oplus^{col} [YA]_2$ where $Z, Y$ are closed
subspaces of $X$, with $Z$ 
a type 2 invariant subspace,
and $\{y^*x : y, x \in Y \} = Y^*Y \subset L^1({\mathcal D})$.  Then
$X$ is simply right $A$-invariant if and only if $Y \neq
\{0\}$.
\item [(2)]  If $X$ is as in {\rm (1)},
then $[Y {\mathcal D}]_2 = X \ominus [XA_0]_2$ (and
$X = [XA_0]_2 \oplus [Y {\mathcal D}]_2$).
\item [(3)]   If $X$ is as described in {\rm (1)},
then that description also holds if $Y$ is replaced by $[Y {\mathcal D}]_2$.  Thus
(after making this replacement)
we may assume that $Y$ is a ${\mathcal D}$-submodule of $X$.
\item [(4)]   The subspaces $[Y {\mathcal D}]_2$ and $Z$ in the decomposition
in  {\rm (1)} are uniquely determined by $X$.  So is $Y$ if we
take it to be a ${\mathcal D}$-submodule (see {\rm (3)}).
\item [(5)]  If $A$ is maximal subdiagonal, then any right $A$-invariant subspace 
$X$ of
$L^2(M)$ is of the form described in {\rm (1)},
with $Y$ the right wandering subspace of $X$.
\end{itemize}
\end{theorem}

\begin{proof}
(1) \ Let $X$ be of the form described above.
We show that $Y \perp
[XA_0]_2$ from which it follows that $X$ is simply right
$A$-invariant if $Y \neq \{0\}$. To see this it is enough to show
that $y \perp (z + xa)$ for any $z \in Z$, $a \in A_0$, and $x,y
\in Y$. From the hypotheses $Y^*Z = \{0\}$ and $Y^*Y \subset
L^1({\mathcal D})$, it now easily follows that $$\tau(y^*(z + xa)) =
\tau((y^*x)a) = \tau(\Phi((y^*x)a)) =  \tau((y^*x)\Phi(a)) = 0$$
as required.
   The converse is obvious.

(3) \
Since $Y^*Z = \{0\}$ and $Y^*Y \subset L^1({\mathcal D})$, we have
${\mathcal D}^*Y^*Z = \{0\}$ and ${\mathcal D}^*Y^*Y {\mathcal D} \subset L^1({\mathcal D})$.
Hence if $\widetilde{Y} = [Y {\mathcal D}]_2$ then
$\widetilde{Y}Z = \{0\}$ and $\widetilde{Y}^*\widetilde{Y}
\subset L^1({\mathcal D})$. Finally, notice that since ${\mathcal D} A = A$ and
$ {\mathcal D} A_0 = A_0$, it is easy to see that $[\widetilde{Y}A]_2 =
[[Y {\mathcal D}]_2A]_2 = [YA]_2$ and $[\widetilde{Y}A_0]_2 = [[Y {\mathcal D}]_2A_0]_2 =
[YA_0]_2$.

(2) \
We saw in the proof of (1) that $Y \perp [XA_0]_2$.
Since $[XA_0]_2 = [XA_0]_2\mathcal{D}$,
we have $[Y {\mathcal D}]_2 \perp [XA_0]_2$. Since $A =
{\mathcal D} + A_0$, it therefore follows that $[YA]_2 = [Y {\mathcal D}]_2 \oplus
[YA_0]_2$. Now since $Z = [ZA_0]_2$, and
$$[YA_0]_2 \subset [XA_0]_2 =
[(Z \oplus [YA]_2)A_0]_2 \subset Z \oplus [YA_0]_2,$$it is clear
that $[XA_0]_2 = Z \oplus [YA_0]_2$. The result therefore follows, since
$$X = Z \oplus
[YA]_2 = Z \oplus ([YA_0]_2 \oplus [Y {\mathcal D}]_2) = [XA_0]_2 \oplus
[Y {\mathcal D}]_2 .$$

(4) \ This follows from (2), and the fact that $Z = X \ominus [YA]_2$.

(5) \  Assume that $A$ is maximal subdiagonal and that $X$ is right
$A$-invariant.  Set $Y = X \ominus [XA_0]_2$. The
subspace $Y$ will clearly be nontrivial if $X$ is simply right
$A$-invariant. We show that $Y^*Y \subset L^1({\mathcal D})$. Let $y, x
\in Y$ be given. Since $[xA_0]_2 \subset [XA_0]_2$ we clearly have
$y \perp [xA_0]_2$, and hence that  $$\tau(y^*xa) = \tau(y^*x(a -
\Phi(a))) + \tau(y^*x\Phi(a)) = \tau(y^*x\Phi(a)) =
\tau(\Phi(y^*x)a)$$for all $a \in A$.   In the last line we have used several properties of $\Phi$ which
are obvious for $\Phi$ considered as a map on $M$, and
which are easily verified for the extension of $\Phi$ to
$L^1(M)$.
(See for example \cite[3.10]{MW}.)
On swapping the roles of $x$ and $y$ and noting that the extension
of $\Phi$ preserves adjoints on $L^1(M)$, we get
$$\tau(y^*xa^*) = \overline{\tau(ax^*y)} =
\overline{\tau(a\Phi(x^*y))} = \tau(\Phi(y^*x)a^*)$$ for all $a
\in A$. So $y^*x - \Phi(y^*x) \perp A + A^*$ which
forces $y^*x = \Phi(y^*x) \in L^1({\mathcal D})$.

Next let $Z = X \ominus [YA]_2$ and let $y \in Y$ and $z \in Z$ be
given. Now by construction $z \perp [yA]_2$ and $y \perp
[zA_0]_2$, whence $$0 = \tau((ya)^*z) = \tau((y^*z)a^*) \quad
\mbox{ for all } a \in A$$ and $$0 = \tau(y^*(za)) = \tau((y^*z)a)
\quad \mbox{ for all } a \in A_0.$$ So $y^*z \perp A + (A_0)^* = A
+ A^*$, which forces $y^*z = 0$.

To see that $ZA \subset Z$, notice that for any $z \in Z$, $y \in
Y$ and $a,b \in A$ it follows from what we just proved that
$(yb)^*(za) = 0$ and hence that $$za \in ([XA]_2 \ominus [YA]_2)
\subset (X \ominus [YA]_2) = Z.$$

Finally, let $V = Z \ominus [ZA_0]_2$. Then since $V \subset Z$, we
have $Y^*V = \{0\}$, which ensures that $V \perp [YA_0]_2$. But by
construction $V \perp [ZA_0]_2$, and so $$V \perp [ZA_0 + YA_0]_2
\supset [(Z \oplus [YA]_2)A_0]_2 = [XA_0]_2.$$But then $V \subset Y = X
\ominus [XA_0]_2$ which by what we noted earlier, forces $V^*V =
\{0\}.$ Clearly $V = \{0\}$, or in other words $Z = [ZA_0]_2$.
\end{proof}

{\bf Remarks.}
1) \ It is an easy consequence of the above Theorem that if $A$ is
maximal subdiagonal and $X$ a right $A$-invariant subspace of $[A]_2$, then
in addition to the other conclusions of Theorem \ref{inv1}
we have that $Z$ is orthogonal to $\mathcal{D}$.  Equivalently, $Z \subset
[A_0]_2$. To see this notice that if $X \subset [A]_2$, then surely $Z \subset
[A]_2$. But then $ZA_0 \subset [A]A_0 \subset [A_0]_2$. Taking the $L^2$-closure
yields $Z = [ZA_0]_2 \subset [A_0]_2$ as claimed. This simple refinement of
the above result is a noncommutative generalization of
Beurling's classical invariant subspace theorem for $H^2$.

2) \ Let $X$ be as in (1) of Theorem \ref{inv1}. Then $Z = [ZA_0]_2 \subset
[XA_0]_2$. Proceeding inductively we conclude that $Z \subset \cap_{n \geq 1}
[XA_0^n]_2$.   From this last fact it follows, for example, that $Z = (0)$ 
if $A$ is the upper triangular matrices.
  
Indeed, for certain maximal subdiagonal algebras, any type 2
invariant subspace is automatically $(0)$, and thus
every closed $A$-invariant subspace $X$ 
is type 1.  Probably one could
isolate various classes of algebras with `shift-like structure' for which
this is also true; in which $n$-fold products
of terms in $A_0$ `converge to zero' in some sense.
(See also \cite{MMS,Zs}).

3) \ For maximal subdiagonal algebras, it is proved in \cite{N} that an
invariant subspace $X$ is of type 1 if and only if
$[X M]_2 = [Y M]_2$, where $Y = X \ominus [X A_0]_2$.
The same result holds for spaces of the form in
Theorem \ref{inv1} (1).

4) \ It is interesting to note that
if $A$ is maximal subdiagonal, then given $f, g \in L^2(M)$ we have that
$f^*g = 0$ if and only if $[fA]_2 \perp [gA]_2$. The ``only if'' part is
obvious, and hence suppose that $[fA]_2 \perp [gA]_2$. It is an
easy exercise to show that this forces $f^*g \perp A+A^*$. Since $A$ is maximal
subdiagonal, this yields $f^*g = 0$.

\bigskip

\begin{proposition} \label{newpr}  Suppose that $X$ is 
as in Theorem {\rm
\ref{inv1}}, and that $W$ is the right wandering subspace of $X$.
Then $W$ may be decomposed as
an orthogonal direct sum $\oplus^2_i \, u_i [{\mathcal D}]_2$,
where $u_i$ are partial isometries in $W \cap M$ with
$u_i^* u_i \in {\mathcal D}$, and $u_j^* u_i = 0$ if
$i \neq j$.   If $W$ has a cyclic vector for the ${\mathcal D}$-action,
then we need only one partial isometry in the above.
\end{proposition}

\begin{proof}   By the theory of representations of a
von Neumann algebra (see e.g.\ the discussion at the
start of Section 3 in
\cite{JS}), any normal Hilbert
${\mathcal D}$-module is an $L^2$ direct sum of cyclic
Hilbert ${\mathcal D}$-modules,
and if $K$ is a normal
cyclic Hilbert ${\mathcal D}$-module, then
$K$ is spatially isomorphic to $[e {\mathcal D}]_2$, for 
an orthogonal projection $e \in {\mathcal D}$. 
 Suppose that the latter isomorphism is implemented by a unitary
${\mathcal D}$-module map $\varphi$.
If in addition $K \subset W$, let $\varphi(e) = 
u \in W$.  Then $\tau(d^* u^* u d) =
\Vert \varphi(e d) \Vert_2^2 = \tau(d^* e d)$, for each $d \in {\mathcal D}$.
By Theorem \ref{inv1}, $u^* u \in L^1({\mathcal D})$,  and so $u^* u = e$.
Hence $u$ is a partial isometry.  Note that $u [{\mathcal D}]_2 \subset
[u {\mathcal D}]_2$ clearly.  However,
$u {\mathcal D} \subset u [{\mathcal D}]_2$, and the latter space is easily seen to
be closed, so that $[u {\mathcal D}]_2 \subset u [{\mathcal D}]_2$.  Thus
$[u {\mathcal D}]_2 = u [{\mathcal D}]_2$.    Putting these facts together,
we see that $W$ is of the desired form.   Note that $u_j^* u_i = 0$ if
$i \neq j$, since  $u_j^* u_i \in {\mathcal D}$,
but $\tau(u_j^* u_i d) = 0$ for any $d \in {\mathcal D}$.
 \end{proof}  

\begin{corollary} \label{adcor}  If $X$ is an
invariant subspace of the form described in
Theorem {\rm \ref{inv1}}, then $X$ is type 1 if and only if
 $X = \oplus^{col}_i \, u_i [A]_2$, for $u_i$ as
in Proposition {\rm \ref{newpr}}.  
\end{corollary}

\begin{proof}  If $X$ is type 1, then $X = [W A]_2$
where $W$ is the right wandering space, and so the one assertion
follows from Proposition {\rm \ref{newpr}}.  
If $X = \oplus^{col}_i \, u_i [A]_2$, for $u_i$ as above, then 
$[X A_0]_2 = \oplus^{col}_i \, u_i [A_0]_2$, and from this
it is easy to argue that  $W = 
\oplus^{col}_i \, u_i [{\mathcal D}]_2$.
Thus $X = [W A]_2 = \oplus^{col}_i \, u_i [A]_2$.
  \end{proof}

 \begin{proposition} \label{typestuff}   Let $X$ be a
closed $A$-invariant subspace of $L^2(M)$, where $A$ is a 
tracial subalgebra of $M$.  
\begin{itemize}  \item [(1)]
If $X = Z \oplus [Y A]_2$ as in Theorem 
 {\rm \ref{inv1}}, then $Z$ is type 2, and $[Y A]_2$ is type 1. 
 \item [(2)]  If $A$ is maximal subdiagonal algebra,
and if $X = K_1 \oplus^{col} K_{2}$ where 
$K_1$ and $K_{2}$ are types 1 and 2 respectively,
then $K_1$ and $K_2$ are respectively the unique
spaces $Z$ and $[Y A]_2$ in  Theorem  {\rm \ref{inv1}}. 
 \item [(3)]  If $A$ and $X$ are as in {\rm (2)}, 
and if $X$ is type 1 (resp.\ type 2), 
then the space  $Z$ of Theorem  {\rm \ref{inv1}}
for $X$ is $(0)$ (resp.\ $Z = X$). 
 \item [(4)]   If  $X = K_1 \oplus^{col} K_{2}$ where
$K_1$ and $K_{2}$ are types 1 and 2 respectively,
then the right wandering subspace for $X$
equals the right wandering subspace for $K_1$.
  \end{itemize} \end{proposition}  

\begin{proof}  (3) \ This is obvious from Theorem {\rm \ref{inv1}}.

(4) \ If $K = K_1 \oplus^{col} K_{2}$ as above,
then $K_{2} = [K_{2} A_0]_2 \subset [K A_0]_2$,
and so $K \ominus [K A_0]_2 \subset K \ominus K_{2} = K_1$.
Thus $K \ominus [K A_0]_2 \subset K_1 \ominus [K_1 A_0]_2$.
Conversely, if $\eta \in K_1 \ominus [K_1 A_0]_2$,
then $\eta  \perp K A_0$ since $\eta^* K_2 = (0)$.
So $\eta \in K \ominus [K A_0]_2$.
  
(1) \  Clearly in this case 
$Z$ is type 2.  To see that $[Y A]_2$ is type 1, note that
since $Y \perp X A_0$, we must have $Y \perp Y A_0$.
Thus $Y \subset [Y A]_2 \ominus [Y A_0]_2$, 
and consequently $[Y A]_2 = [([Y A]_2 \ominus [Y A_0]_2) A]_2$.

(2) \  Suppose that $X = K_1 \oplus^{col} K_{2}$ where
$K_1$ and $K_{2}$ are types 1 and 2 respectively.
Let $Y$ be the right wandering space for $K_1$.
By Theorem \ref{inv1} we have $Y^* Y \subset L^1({\mathcal D})$.
So $X = [Y A]_2 \oplus^{col}  K_{2}$,
and by the uniqueness assertion in Theorem \ref{inv1}, 
$K_2$ is the space $Z$ in Theorem \ref{inv1} for $X$,
and $K_1 = K_{2}^\perp = [Y A]_2$.
\end{proof}

In the case $p = 2$, items (1)--(3) in Theorem \ref{main}
 follow from the last results.

\bigskip
  
{\bf {\em Proof} of Corollary \ref{stand} if $p = 2$:}
  If $K = u H^2$ for a unitary $u$, then as in 
the proof of Corollary \ref{adcor},
$K \ominus [K A_0]_2 = u [{\mathcal D}]_2,$
which has separating and cyclic vector $u$.  Conversely,
if the right wandering subspace $W$ has a
separating and cyclic vector $v$ for ${\mathcal D}$,
then  the  proof of Proposition \ref{newpr}, with $e = 1$ in
that proof, shows that $u^* u = 1$.  So $u$ is unitary.
Since $K_2^* u = (0)$ we have $K_2  = (0)$.  Thus
$K = [W A]_2 = u H^2$.
 \hfill $\square$

\bigskip

The equivalences in the next result generalize, and 
give as an immediate consequence, the equivalence of (a)
and (c) in \ref{rest}, in the antisymmetric case.

\begin{corollary}  \label{Co}  
 Let $A$ be a tracial algebra.  The following are equivalent:
\begin{itemize} \item [(i)]  $A$ is maximal subdiagonal,
  \item [(ii)]  For every right
$A$-invariant subspace $X$ of $L^2(M)$, the right wandering subspace 
 $W$ of $X$ satisfies
$W^* W \subset L^1({\mathcal D})$, and $W^* (X \ominus [W A]_2) = (0)$.
\item [(iii)]  Every right
invariant subspace of $L^2(M)$ satisfies {\rm (1)} and  {\rm (3)}
of  Theorem {\rm \ref{main}}.
 \end{itemize} \end{corollary}

\begin{proof}  The fact that (i) implies (ii) is
proved in Theorem \ref{inv1}.  The fact that (ii) implies (i) may be proved via
the later Beurling-Nevanlinna type factorization theorem
\ref{NB}, along the same lines as the proof of Proposition \ref{rest}.
We choose to also give a direct proof.
To see that  (ii) implies (i), first set $X = L^2
\ominus [A_0^*]_2$. We will deduce that $A$ satisfies
$L^2$-density. That is, $L^2(M)$ is the closure of $A + A^*
= A + A_0^*$, or equivalently that $X = [A]_2$. 
To this end, note that $X$ is right A-invariant. 
It is easy to see that $1 \in W$, which forces $X \ominus [W A]_2 = (0)$,
 and $W \subset W^* W  \subset L^1({\mathcal D})$.
Thus $W \subset L^2(M) \cap L^1({\mathcal D}) = L^2({\mathcal D})$. So $X =
[WA]_2 \subset [A]_2$. The converse inclusion $[A]_2 \subset X$
follows from the fact that $[A]_2$ is orthogonal to $[A^*_0]_2$.

We now prove that $A$ possesses the unique normal state extension
property; so that
 $A$ is maximal subdiagonal.  To this end, let $g \in
L^1(M)_+$ with $\tau(g A_0) = 0$.
 We may assume that $g \neq 0$. Let $h = g^{\frac{1}{2}}
\in L^2(M)$, and set $X = [h A]_2$. 
Note that $h \perp [h A_0]_2$ since if $a_n
\in A_0$ with $h a_n \to k$ in $L^2$-norm, then $\tau(h^* k) =
\lim_n \tau(h^* h a_n) = 0$. In particular, the fact that $h \perp
[h A_0]_2$ ensures that $h \in X \ominus [XA_0]_2 = W$.  By
hypothesis, $h^* h = g \in L^1({\mathcal D})$.

We have seen already that (i) implies (iii).
If (iii) holds then by Proposition \ref{typestuff} (4) 
the wandering subspace $W \subset K_1$, and $K \ominus [W A]_2 = 
K \ominus K_1 = K_2$.  Thus  $W^* K_2 \subset
K_1^* K_2 = (0)$, and $W^* W \subset L^1({\mathcal D})$.      
This is (ii).
\end{proof}

{\bf Remark.}   The conditions in the last result are also equivalent to
(iv) \ every simply right
$A$-invariant subspace $X$ of $L^2(M)$ is of the form $X = Z
\oplus^{col} [YA]_2$ where $Y, Z$ are closed subspaces of $X$ with $Z$
type 2, and $(0) \neq Y^*Y \subset L^1({\mathcal D})$;
and to (v) \ every right
$A$-invariant subspace of $L^2(M)$ satisfies {\rm (1)} and  {\rm (2)}
of  Theorem {\rm \ref{main}}.  This is fairly obvious from the
proof of \ref{Co}.
  
\section{Noncommutative Beurling-Nevanlinna factorization}

We now discuss several general inner-outer,
or Beurling-Nevanlinna, factorization theorems
each of which turns out to characterize maximal subdiagonal algebras.
One of these theorems (Theorem \ref{NB} (iv)),
has the classical equivalence with the Beurling-Nevanlinna
factorization theorem as an immediate consequence, or special case
(namely that (b) is equivalent to (c) in Proposition
\ref{rest} in the classical situation).  This will be clear from 
an argument in the present paragraph.
Our other Beurling-Nevanlinna factorization theorems
 (namely Theorem \ref{NB} (ii) or (iii)) have the advantage of
having more attractive hypotheses, which are perhaps easier to
check.  To understand the hypothesis of these theorems, suppose
for a moment that ${\mathcal D}$ is one dimensional
as in the classical case,
and that $f \in L^2(M)$ is such that $f \notin [f A_0]_2$.
Thus  $[f
A]_2$ is type 1.
The right wandering subspace $[f
A]_2 \ominus [f A_0]_2$ in this case is
also one dimensional (since $f A \subset [f A_0]_2 + f
{\mathcal D} \subset [f A]_2$, and since
$[f A_0]_2 + f
{\mathcal D}$ is closed,  $[f A_0]_2$ has codimension
one in $[f A]_2$).  In particular the wandering subspace has a
separating and cyclic vector.  In the general case,
having such a vector is a necessary condition for BN-factorizability:

\begin{lemma}  \label{res}    Suppose that $A$ is a tracial
algebra, and that $f \in L^2(M)$ is BN-factorizable (resp.\
partially BN-factorizable).
Then the right wandering subspace of $[f A]_2$ has
a separating and cyclic vector (resp.\ a
cyclic vector) for the right ${\mathcal D}$-action.  
In fact, that vector may be taken to be a unitary
(resp.\ partial isometry) in $M$.
Conversely, if $A$ is maximal subdiagonal,
and if the right wandering subspace of $[f A]_2$ has
a separating and cyclic vector,
then $f$ is BN-factorizable.
\end{lemma}

\begin{proof}  If $f = u h$ is a BN-factorization, then
$[fA]_2 = [u h A]_2 = u [h A]_2 = u [A]_2$.   Similarly,
$[fA_0]_2 = u [h A_0]_2 = u [A_0]_2$, the latter since
$h A_0  \subset [A]_2 A_0 \subset [A_0]_2$ and
$$A_0 = A A_0 \subset [A]_2 A_0 \subset [h A]_2 A_0 \subset [h A_0]_2 .$$
Thus $[f A]_2
 \ominus [f A_0]_2 = u ([A]_2 \ominus [A_0]_2) = u [{\mathcal D}]_2$.
Similar, but slightly more cumbersome,  arguments work in
the partial isometry case.

Suppose that $A$ is maximal subdiagonal and the right wandering subspace 
of $[f A]_2$ has a separating and cyclic vector,
then by Corollary \ref{stand},
$[f A]_2 = u [A]_2$ for a unitary $u$.
We may write $f = u h$ for $h \in [A]_2$.
Clearly $u \in M \cap [f A]_2$, and 
$$[h A]_2 = u^* u [h A]_2
= u^* [f A]_2
 = [A]_2 .$$
This completes the proof.
\end{proof}

Recall that a {\em wandering vector} is a vector $f \in L^2(M)$ with
$f \perp [f A_0]_2$.  Examples include the partial isometries $u_i$ in 
the previous section.
Wandering vectors were completely characterized in
\cite{AIOA} in the case that $A$ is maximal subdiagonal.
Note that if $f$ is a wandering vector,
 then $[f {\mathcal D}]_2 \perp [f A_0]_2$, and so one easily sees that
$[f {\mathcal D}]_2$ is the right wandering subspace of $[f A]_2$.
Moreover, $[([f A]_2 \ominus [f A_0]_2) A]_2 =
[[f {\mathcal D}]_2 A]_2 = [f A]_2$, so that $[f A]_2$ is type
1.   In this
case, to say that the right wandering subspace
of $[f A]_2$ has a cyclic separating
vector is equivalent to saying that $[f {\mathcal D}]_2$
has such a vector.  By \cite[Exercise 9.6.2]{KR}, this
is also equivalent to saying that $f$ is separating; that is,   the
map $d \mapsto f d$ on ${\mathcal D}$ is one-to-one.

\begin{proposition} \label{BN}  Suppose that $A$ is a tracial
algebra.
If $0 \neq f \in L^2(M)$,
consider the following conditions:
\begin{itemize}
\item [(a)]  $f \in M^{-1} \cap M_+$,
\item [(b)]  $f$ is a wandering vector, and the
map $d \mapsto f d$ on ${\mathcal D}$ is one-to-one,
\item [(c)]  the right wandering subspace
of $[f A]_2$ has
a nonzero separating and cyclic vector for the right action of ${\mathcal D}$,
\end{itemize}
Then any one of conditions {\rm (a)} or {\rm (b)}
imply {\rm (c)}.
\end{proposition}

\begin{proof}  
Suppose that {\rm (b)} holds.  In this case,
$[f {\mathcal D}]_2$ is the right wandering subspace
of $[f A]_2$, as remarked above, and so {\rm (c)} holds.

If {\rm (a)} holds, indeed if $f \in M^{-1}$,
then $[f A]_2 = f [A]_2$ and $[f A_0]_2 = f [A_0]_2$.
Thus as purely algebraic ${\mathcal D}$-modules,
$$[f A]_2 \ominus [f A_0]_2 \cong [f A]_2 / [f A_0]_2  \cong
[A]_2 / [A_0]_2 \cong [{\mathcal D}]_2 .$$
The latter module has a nonzero separating and cyclic vector.
Since these isomorphisms are also continuous,  so does
$[f A]_2 \ominus [f A_0]_2$.  So {\rm (c)} holds.  
\end{proof}

\begin{theorem} \label{NB} \ (Beurling-Nevanlinna factorization for
tracial algebras) \
For a tracial subalgebra $A$ of $M$, the following are equivalent:
\begin{itemize}
\item [(i)]  $A$ is maximal subdiagonal,
\item [(ii)]  Every nonzero $f \in L^2(M)$ satisfying either
{\rm (a)} or {\rm (b)}
in  {\rm \ref{BN}} is BN-factorizable,
\item [(iii)]  If $f \in M^{-1} \cap M_+$, or if $f$ is a wandering vector,
then $f$ is partially BN-factorizable,
\item [(iv)]   Every $f \in L^2(M)$ satisfying {\rm (c)} in  {\rm \ref{BN}} is BN-factorizable.
\end{itemize}
Indeed, every $f \in M^{-1} \cap M_+$ is BN-factorizable,
if and only if $A_\infty$ is maximal subdiagonal.
Also, $A$ has the unique normal state extension property if and only if
every wandering vector is partially BN-factorizable.
\end{theorem}

\begin{proof}
The proof of \cite[Lemma 2.4.3]{Sr} demonstrates one direction of the
penultimate `if and only if', and we leave the converse as an exercise.

For the equivalences between (i), (ii), and (iv),
by Proposition \ref{BN} and Lemma \ref{res}
it remains to show that (ii)
implies (i). If (ii) holds, then by the last statement $A_\infty$
is maximal subdiagonal.
Hence $A$ satisfies $L^2$-density. If therefore we can show that
$A$ also satisfies the unique normal state extension property we
are done (as in the proof of
Corollary \ref{Co}). Suppose that $g \in L^1(M)_+$ satisfies $\tau(g A_0) =
0$.
We need to show that $g \in L^1({\mathcal D})_+$. Since $\tau((g + 1) A_0) =
0$, we can replace $g$ with $g + 1$ if necessary.
Let $f = g^{\frac{1}{2}} \in L^2(M)$. Then $f \perp [f A_0]_2$.
Obviously, $f$ is a cyclic vector for $[f {\mathcal D}]_2$. If $f d 
= 0$ then
$d = 0$ since $\tau(d^* d) \leq \tau(d^* f^* f d)$. So $f$ is
separating.
By hypothesis, $f = u h$ for an outer  $h \in [A]_2$ and
some unitary $u$ in $M$.   Since $h = u^* f \perp [h A_0]_2
= [A_0]_2$, and $h \in [A]_2$,
it follows that $h \in [{\mathcal D}]_2$.   Thus $g \in [{\mathcal D}]_1 = L^1({\mathcal D})$.

An easy modification of the last argument gives one direction of the
last assertion of the Theorem.
For the other direction, suppose that $A$ has the unique normal state extension property,
and $f \in L^2(M) \ominus [f A_0]_2$.   By the remarks preceding
Proposition \ref{BN}, $[f A]_2$ is a type 1 invariant subspace.
Following the first half of the proof of Theorem \ref{inv1} (5), we see that
$f^* f \subset L^1({\mathcal D})$, and so $[f A]_2$ has
a decomposition of the type considered in
Theorem \ref{inv1} (here $Z = (0)$).   Since $Y = [f {\mathcal D}]_2$ is
cyclic, by Proposition \ref{newpr} there is a partial isometry $u$ with
$[f A]_2 = u [A]_2$ and $u^* u \in {\mathcal D}$.   Clearly $u \in [f A]_2 \cap M$,
We may write $f = u h$ with $h \in (u^* u) [A]_2 \subset [A]_2$.
Also, $u [h A]_2 = [f A]_2 = u [A]_2$, so that
$[h A]_2 = (u^* u) [h A]_2 = (u^* u) [A]_2$.  Thus $f$ is partially
BN-factorizable.

The equivalence of (i) with (iii) follows from the last two
iff's of our theorem's statement, 
the fact that if $A_\infty$ is maximal subdiagonal
then $A$ satisfies $L^2$-density, and the fact that $f \in M^{-1}$
is  partially BN-factorizable if and only if it is BN-factorizable.
Indeed if $f$ is invertible, with partial BN-factorization
$f = u h$, then $h = u^* f$ is bounded, and
so both $u$ and $h$ are also invertible since we are in a finite von
Neumann algebra.  Hence $u$ is a unitary. 
  \end{proof}

At the end of Section 4 we give
 a very general `inner-outer factorization' result.

\section{The case of $L^p$ for $p \neq 2$}

To discuss the $L^p$-version of some of the results above,
we use the  `column $L^p$-sum'
from \cite{JS}.  
Suppose that $1 \leq p < \infty$ and that 
 $\{ X_i : i \in I \}$ is
a collection of closed subspaces of $L^p(M)$.  We then define
the {\em external} column $L^p$-sum $\oplus^{col}_i \, X_i$ to be the
closure of the restricted algebraic sum in the norm
$\Vert (x_i) \Vert_p \overset{def}{=} \tau((\sum_i \, x_i^* x_i)^{\frac{p}{2}})^{\frac{1}{p}}$.
That this is a norm for $1 \leq p < \infty$ is verified in
\cite{JS}.   If $X$ is a subspace of $L^p(M)$, and if  $\{ X_i : i \in I \}$ is
a collection of  subspaces of $X$, which together densely
span $X$, with the
property that $X_i^* X_j = \{ 0 \}$ if $i \neq j$, then
we say that $X$ is the {\em internal} column $L^p$-sum $\oplus^{col}_i \, X_i$.
Note that  if $J$ is a finite subset of $I$,
and if $x_i \in X_i$ for all $i \in J$, then we have that
$$\tau(|\sum_{i \in J} \, x_i|^p)^{\frac{1}{p}}
= \tau((|\sum_{i \in J} \, x_i|^2)^{\frac{p}{2}})^{\frac{1}{p}}
= \tau((\sum_{i \in J} \, x_i^* x_i)^{\frac{p}{2}})^{\frac{1}{p}} .$$
This shows that $X$ is then isometrically isomorphic
 to the external column $L^p$-sum $\oplus^{col}_i \, X_i$.
Since the projections onto the summands are clearly
contractive, it follows by routine arguments
(or by \cite[Lemma 2.4]{JS}) that if $(x_i) 
\in \oplus^{col}_i \, X_i$, then the net $(\sum_{j \in J} \, x_j)$,
indexed by the finite subsets $J$ of $I$, converges in norm to $(x_i)$.

We will need a couple of technical results:

\begin{lemma}   \label{istr}
If $0 < p < \infty$, and if $e$ is a projection in $M$,
 then $v \in L^p(M)_+$ satisfies
$\tau((d^* v d)^p) = \tau((d^* e d)^{p})$ for all $d \in M$,
if and only if $v = e$.
\end{lemma}

\begin{proof}   We regard $v$ as an unbounded
positive operator affiliated
with $M$, then $v^p \in L^1(M)_+$.  Choosing $d = e^\perp$ shows that
$e^\perp v e^\perp = 0$.  It is easy to see, as in the bounded
case, that these imply that $v e^\perp = e^\perp v = 0$ too.
So $v = e v e$.   Replacing $M$ by $e M e$, we can assume that $e = 1$.
If $E$ is a projection in the minimal
von Neumann algebra $M_0$ generated by $v$ (see e.g.\ \cite[p.\ 349]{KR}),
then $(E v E)^p =
E v^p$, and so by hypothesis
$\tau(E (v^p - 1)) = 0$.  If $E$ is the spectral projection 
for $v$ corresponding
to $[0,1]$, then $E v^p \leq E$ and so we have $E = E v^p$ since $\tau$ is
faithful.  On the other hand, $E^\perp v^p \geq E^\perp$ so that $E^\perp
 = E^\perp v^p$.  Thus $v^p = 1$.
\end{proof}

\begin{lemma} \label{ll1}   Let $A$ be maximal subdiagonal and let $K$
be an $A$-invariant subspace of $L^p(M)$, 
for $2 \leq p \leq \infty$.
  Then $[K]_1 \cap L^p(M) = [K]_2 \cap L^p(M) = [K]_p$.
(The last symbol, if $p = \infty$, is always taken in this paper
to mean the weak* closure.)
\end{lemma}

\begin{proof}  In this proof we will denote polars taken with respect to
the dual pair $(L^1(M),M)$ by $\perp$,
and polars taken with respect to
the dual pair $(L^p(M),L^q(M))$ by $\circ$.
Here $\frac{1}{p} + \frac{1}{q} = 1$.  We assume that
$p < \infty$, and leave the case $p = \infty$ to the reader.
  Note that if $K$ is a right $A$-invariant subspace of $L^p(M)$ then 
$K^\circ$ is a closed right $A^*$-invariant subspace of $L^{q}(M)$. To see this
note that we clearly have
$$0 = \tau(y^*(xa)) = \tau((ya^*)^*x) , \qquad x \in K, a \in A, y \in K^\circ,$$
since  $xa \in K$. But then $ya^* \in K^\circ$ as
required.

Since $K^\circ$ is a closed $A^*$-invariant subspace of $L^{q}(M)$,
by the result of Saito mentioned in our introduction
 $K^\circ \cap M$ is norm dense in $K^\circ$.
Regarding $K$  as a subspace of $L^1(M)$ we have
$K^\perp = K^\circ \cap M$. Hence by the bipolar
theorem $(K^\circ \cap M)_\perp = [K]_1$.
Clearly $(K^\circ \cap M)_\perp \cap L^p(M) = (K^\circ \cap M)^\circ$.
Hence
 $$[K]_1 \cap L^p(M)
 = (K^\circ \cap M)_\perp \cap L^p(M) = (K^\circ \cap M)^\circ = 
(K^\circ)^\circ = [K]_p.$$
The other assertion now follows from the fact that
$K \subset [K]_2 \cap L^p(M) \subset [K]_1\cap L^p(M)$.
If $p = \infty$ it is easy to check that 
$[K]_2 \cap M$ is weak* closed.  So in all cases,
$[K]_p  \subset [K]_2 \cap L^p(M)$, which gives the result.
   \end{proof}

\begin{corollary} \label{bij}  If $A$ is maximal subdiagonal
then for any $1 \leq p \leq q \leq \infty$ there is a lattice isomorphism
between the closed (weak*-closed, if $q = \infty$)
 right $A$-invariant subspaces of $L^p(M)$ and $L^q(M)$.
\end{corollary}

\begin{proof}  We may take $p$ or $q$ to be 2.  The isomorphism (resp.\ its inverse) of 
course is the map taking $K$ to its closure (resp.\ intersection with 
the appropriate $L^p$ space).  
This follows easily from the aforementioned result of Saito,
and the Lemma (or a tiny variant of it).
\end{proof}
 
\begin{definition}  We define the right wandering subspace of $K$,
if $1 \leq p \leq 2$ (resp.\
$p \geq 2$) to be the $L^p$-closure of the 
right wandering subspace of $K \cap L^2(M)$
(resp.\ to be the intersection of $L^p(M)$ with the
right wandering subspace of $[K]_2$).
\end{definition}
  
\begin{theorem} \label{lp}  Let $A$ be a maximal subdiagonal subalgebra
of $M$, and suppose that $K$ is a closed $A$-invariant 
subspace of $L^p(M)$,
for $1 \leq p \leq \infty$.  (For $p = \infty$ we assume that
$K$ is weak* closed.)
\begin{itemize}  \item [(1)]  $K$ may be written as a column $L^p$-sum
$K = Z \oplus^{col} (\oplus_i^{col} \, u_i H^p)$, where
$Z$ is a closed (indeed weak* closed if $p = \infty$)
type 2 subspace of $L^p(M)$,
and where $u_i$ are partial isometries in $M \cap K$ with $u^*_j
u_i = 0$ if $i \neq j$,
and with $u_i^* u_i \in {\mathcal D}$.
Moreover, for each $i$, $u_i^* Z = (0)$, left multiplication
by the $u_i u_i^*$ are contractive projections from $K$ onto
the summands $u_i [A]_p$, and left multiplication
by $1 - \sum_i \, u_i u_i^*$ is a contractive projection from $K$ onto $Z$.
\item [(2)]   The wandering quotient $K/[K A_0]_p$ is isometrically
${\mathcal D}$-isomorphic to the right wandering subspace of $K$;
and the latter also equals 
$\oplus_i^{col} \, u_i [{\mathcal D}]_p$, where 
$u_i$ are from {\rm (1)}.  
(Here $[\cdot]_\infty$ is the weak* closure as usual.)
\item [(3)]   $K$ is type 1 if and only if
 $K  \cap L^2(M)$ (resp.\ $[K]_2$) is type 1,
and if and only if $Z = (0)$ in {\rm (1)}.
If $1 \leq p \leq 2$ (resp.\
$p \geq 2$), then $K$ is type 2 if and only if
$K  \cap L^2(M)$ (resp.\ $[K]_2$) is type 2, and if and only if $K = Z$.
\end{itemize} 
\end{theorem}

\begin{proof} 
(1) \  First suppose that $p \leq 2$.
By Saito's theorem  mentioned in the introduction, $K$ is the
$L^p$-closure of $K  \cap L^2(M)$.  Theorem
\ref{inv1} gives a decomposition $K  \cap L^2(M) = Z' \oplus [Y A]_2$,
where $Z'$ is type 2, and $Y$ is the right 
wandering subspace of $K  \cap L^2(M)$.
 Let  $Z$ be the $L^p$-closure of $Z'$.
Note that  $Z$ is type 2.  Indeed $[Z A_0]_p = [[Z']_p A_0]_p
= [Z' A_0]_p$, as is easy to check.  Thus
$$[Z A_0]_p =  [Z' A_0]_p = [[Z' A_0]_2]_p = [Z']_p = Z.$$
We leave it as an exercise that
$Z \cap L^2(M) = Z'$. 

Since $z^* y = 0$ for $y \in Z'$, $y \in Y A$, it follows
using  (4.8) in \cite{JS}
 that $(0) = Z^* [Y A]_p$.
From this it is easy to see that $K = Z \oplus^{col} K_1$
where $K_1 = [Y A]_p$.
By Corollary \ref{adcor}, $[Y A]_2 =
\oplus_i^{col} \;  u_i [A]_2$, for $u_i$ as above.
Thus the  $L^p$-closure of $\sum_i \;  u_i A$ is
all of $K_1$.  On the other hand, $(u_j [A]_p)^* (u_i [A]_p) =
(0)$ if $i \neq j$.   So $K_1 = \oplus^{col}_i \, u_i [A]_p$.

Now suppose that $2 \leq p \leq \infty$.
 We embed $K$ into $L^2(M)$; by the $L^2$ result
$[K]_2 = Z' \oplus^{col} (\oplus_i^{col} \, u_i H^2)$ as usual,
where $Z'$ is type 2.   We assume that $p < \infty$ in
the arguments below, however the case $p = \infty$ is
a simple variant.    Let $K_1 =
(\oplus_i^{col} \, u_i H^2) \cap L^p(M)$.   
Note that $[\oplus_i^{col} \, u_i H^p]_2
= \oplus_i^{col} \, u_i H^2$,
so that by Lemma \ref{ll1} we have $K_1 = \oplus_i^{col} \, u_i H^p$.

Define $Z = Z' \cap L^p(M)$, which is also 
easily seen to be closed.   Since $Z'$ 
is $A$-invariant,  by Saito's result $Z$ 
is an $A$-invariant subspace of $L^p(M)$ which is $L^2$-dense in $Z'$.
  Set $X_0 = Z \oplus K_1$, it follows that
$X_0$ is $L^2$-dense in $[K]_2$.  Also, $X_0$ is
an invariant subspace of $L^p(M)$, which is easily seen 
to be closed.
Clearly $Z^* K_1 = (0)$.  By Lemma \ref{ll1} we 
have
$$K  = [K]_2 \cap L^p(M) = [X_0]_2 \cap L^p(M) = X_0 = Z \oplus^{col} K_1.$$
To see that $Z$ is type 2,  by 
Lemma \ref{ll1}, and the fact that $[Z]_2 = Z'$, we have  
$$Z = Z' \cap L^p(M) = [Z' A_0]_2 \cap L^p(M)
\subset  [Z A_0]_2 \cap L^p(M) = [Z A_0]_p. $$

Finally, suppose $1 \leq p \leq \infty$.
Since left multiplication by $u_i u_i^*$
annihilates $Z'$ and $u_j [A]_p$ if $j \neq i$, left multiplication
by the $u_i u_i^*$ are contractive projections from $K$ onto
the summands $u_i [A]_p$, and left multiplication
by $1-\sum_i \, u_i u_i^*$ is a contractive projection onto $Z$.

(2) \   The fact that the right wandering subspace
equals $\oplus^{col}_i \, u_i [{\mathcal D}]_p$ follows
from a slight modification of the argument in (1) 
that $K_1 = \oplus^{col}_i \, u_i H^p$.    We also need the fact that 
$L^p(M) \cap [{\mathcal D}]_1 = [{\mathcal D}]_p$.
  
By e.g.\ \cite{MW}, $\Phi$ induces a 
contractive `expectation' on every $L^p(M)$. 
Assume $p \geq 2$.  Then $\Phi(x)^* \Phi(x) \leq \Phi(x^* x)$ for $x \in L^p(M)$,
as may be seen by routine continuity arguments.
Again we assume $p < \infty$, and leave the 
variation of the argument in the case $p = \infty$ to the 
reader.   Define a map $\theta : K = K_1 \oplus^{col} K_2 \to
K$ by
$\theta(w) = \sum_i \, u_i \Phi(x_i)$
if $w = \sum_i \, u_i x_i + k_2$ for $x_i \in
H^p$ and $k_2 \in K_2$.   It is easy to see that
$\theta(w)$
equals $\sum_i \, u_i \Phi(u_i^* w)$, which shows that
$\theta$ is well defined (and weak* continuous
if $p = \infty$).  Indeed since $u_i^*u_i \in {\mathcal D}$, we have
$$\tau((\sum_i \, \Phi(x_i)^* u_i^* u_i \Phi(x_i))^{\frac{p}{2}})
\leq \tau(\Phi(\sum_i \, x_i^* u_i^* u_i x_i)^{\frac{p}{2}})
\leq \tau((\sum_i \, x_i^* u_i^* u_i x_i)^{\frac{p}{2}}),$$
which shows that $\Vert \theta(w)  \Vert_p^p
\leq \Vert  w \Vert_p^p.$    Thus $\theta$ is
a contractive projection onto its range, and hence
induces an isometric ${\mathcal D}$-module map from
$K/{\rm Ker}(\theta)$ onto $\oplus^{col}_i \, u_i [{\mathcal D}]_p$.
Since $x_i A_0 \in [A_0]_p \subset {\rm Ker}(\Phi)$ if
$x_i \in H^p$, it is clear that ${\rm Ker}(\theta)$ contains $K A_0$,
and hence contains $[K A_0]_p$.  Conversely
if $\theta(\sum_i \, u_i x_i + k_2) = 0$ then
$u_i^* u_i \Phi(x_i) = 0$ for every $i$.
Thus $u_i^* u_i x_i \in {\rm Ker}(\Phi) \cap H^p = [A_0]_p$.
Thus $u_i x_i \in u_i [A_0]_p$, and so
$\sum_i \, u_i x_i + k_2 \in [K A_0]_p$.  So  ${\rm Ker}(\theta)
= [K A_0]_p$.

If $1 \leq p < 2$, we follow the same argument with the 
following modification.  
Suppose that $1/p + 1/q = 1$, and that $F$ is a finite set
of indices.  Define $\theta_F$ on $L^q(M)$ taking
$w \mapsto \sum_{i \in F} \, u_i \Phi(u_i^* w)$.
Since $\sum_F u_i u_i^*$ is a projection, by the arguments in the last
paragraph
it is easy to see that $\tau((\sum_F \, \Phi(u_i^*w)^* u_i^* u_i
\Phi(u_i^*w))^{\frac{q}{2}})$ is dominated by
$$\tau(\Phi(w^* (\sum_F \, u_i
u_i^*) w)^{\frac{q}{2}}) \leq  \tau(\Phi(w^*w)^{\frac{q}{2}}) \leq
\tau((w^*w)^{\frac{q}{2}}) = \tau(|w|^q),$$
which shows that $\theta_F$ is a contraction.
For $g \in L^p(M), f \in L^q(M)$ we have that
$$\tau(g^* \theta_F(f))
= \sum_{i \in F} \, \tau(g^* u_i \Phi(u_i^* f))
= \sum_{i \in F} \, \tau(\Phi(g^* u_i) u_i^* f)
= \tau(( \sum_{i \in F} \, u_i \Phi(u_i^* g))^* f),$$ 
where the middle equality follows by e.g.\ 
\cite[3.10]{MW}.
This shows that the contraction
$(\theta_F)_* \in B(L^p(M))$ 
is  precisely the map $w \mapsto \sum_{i \in F} \, u_i \Phi(u_i^* w)$
on $L^p(M)$.
Since this holds for every finite subset $F$,
this implies that the map $\theta$ above is a densely defined
contraction on $K$, and thus extends continuously to $K$.
By continuity it follows that this extension has precisely the same 
formula as before, and now the earlier argument works.   
 
(3) \  If $K \cap L^2(M)$ (resp.\ $[K]_2$) is type 1 then it is obvious
from the proof of (1) that $Z' = Z = (0)$. 
Similarly, if $K \cap L^2(M)$ (resp.\ $[K]_2$) is type 2
then it is obvious that $K = Z$.
It is trivial that if $K = Z$ then $K$ is type 2.
Conversely, if $K$ is type 2 then the wandering quotient is
$(0)$.  Identifying
this  with the subspace of $K$ described in (2), all the $u_i$
are zero.  Thus $K = Z$.  Also, if $p > 2$ then $[K]_2 = [Z]_2 = Z'$
which is type $2$; and a similar argument works if
$p < 2$.    Since $[(\oplus^{col}_i \, u_i [{\mathcal D}]_p) A]_p
= \oplus^{col}_i \, u_i H^p$, we have that $Z = (0)$
iff $\oplus^{col}_i \, u_i [{\mathcal D}]_p$ generates
$K$.   By (2) this happens iff $Z$ is type 1.
 
If
$1 \leq p \leq 2$,
suppose that $K = \oplus^{col}_i \, u_i H^p$ for partial isometries $u_i$ 
satisfying the relations in (1).  
Then $[\oplus^{col}_i \, u_i H^2]_p = K$, so by 
Lemma \ref{ll1} we have
$K \cap L^2(M) \subset [\oplus^{col}_i \, u_i H^2]_1 \cap L^2(M)
= \oplus^{col}_i \, u_i H^2$.  So
$\oplus^{col}_i \, u_i [A]_2 = K \cap L^2(M)$,
which is type 1 by Corollary \ref{stand}. 

If  $2 \leq p \leq \infty$
and $K = \oplus^{col}_i \, 
u_i H^p$ for partial isometries $u_i$
satisfying the relations in (1), then 
it is easy to argue that
$\oplus^{col}_i \, u_i [A]_2 = [K]_2$.  So $[K]_2$ is
type 1.   (In the case $p = \infty$,
note that $\sum_i \, u_i A \subset [K]_2 \cap M = K$
by Lemma \ref{ll1}, so that $[\overline{\sum_i \, u_i A}^{weak*}]_2
 \subset [K]_2.$)  
\end{proof}

We have now proved the existence of
the type decomposition in Theorem \ref{main} (1).  
The uniqueness of this type decomposition follows
from the following:

\begin{corollary} \label{untip}  If $K$ is a subspace of $L^p(M)$ of the 
form $K = K_1 \oplus^{col} K_2$ where $K_1$ is type 1
and $K_2$ is type 2, and if  $1 \leq p \leq 2$, then
$K_1$ (resp.\ $K_2$) is the $L^p$ closure of the type 1 
(resp.\ type 2) part of
$K \cap L^2(M)$.
If $2 \leq p \leq \infty$ then
$K_1$   (resp.\ $K_2$) is the intersection of $L^p(M)$ 
with the type 1
(resp.\ type 2) part of $[K]_2$.   \end{corollary}

\begin{proof}
If $1 \leq p \leq 2$, then 
clearly $(K_1 \cap L^2(M)) + (K_2 \cap L^2(M)) \subset 
K \cap L^2(M)$.  On the other hand, if 
$x \in K \cap L^2(M)$, we can write $x = k_1 + k_2$,
for $k_i \in K_i$.  Since $K_1$ is type 1,
we can write $K_1 = \oplus^{col}_i \, u_i [A]_p$ for 
some $u_i$ as above.  So $k_1 = 
\sum_i \, u_i u_i^* k_1$.  Since 
$u_i \in K_1$ we have $u_i^* k_2 = 0$.
Thus since $x = (1-\sum_i \, u_i u_i^*) x + \sum_i \, u_i u_i^* x$,
we have $k_1 = \sum_i \, u_i u_i^* x$, and so
$k_2 = (1-\sum_i \, u_i u_i^*) x$.  This forces 
$k_1 \in (\sum_i \, u_i u_i^*) L^2(M) \subset L^2(M)$,
similarly $k_2  \in L^2(M)$.  Thus we have
$K \cap L^2(M) = (K_1 \cap L^2(M)) \oplus^{col} (K_2 \cap L^2(M))$.
Since $K_1 \cap L^2(M)$ is type 1 and $K_2 \cap L^2(M)$ is
type 2, we have by the uniqueness of decomposition in
the $L^2$ case, that
$K_1 \cap L^2(M) = [(K \cap L^2(M) \ominus [(K \cap L^2(M)) A_0]_2) A]_2$.
By Saito's theorem $K_1$ is the $L^p$ closure of
$K_1 \cap L^2(M)$.  Since $K_2 = \{ z \in K  : z^* K_1 = (0) \}$,
it is uniquely determined, and so the 
assertion about $K_2$ follows by the proof of the previous theorem.

If $2  \leq p$, then $K = K_1 \oplus^{col} K_2$
implies easily that $[K]_2 = [K_1]_2 \oplus^{col} [K_2]_2$.
Since $[K_1]_2$ and $[K_2]_2$ are types 1 and 2 respectively
(by (3) of the Theorem), the result follows from Lemma \ref{ll1}.
\end{proof}

\begin{proposition} \label{triple}  Suppose that $K$ is a
closed (indeed weak* closed, if $p = \infty$)
$A$-invariant subspace of $L^p(M)$, and 
that  $A$ is maximal subdiagonal.
The following are equivalent:
\begin{itemize} \item [(1)]  $K = u H^p$ for a unitary $u \in M$.
 \item [(2)]  The right wandering
subspace of $K$ (or equivalently, the right wandering
quotient) has a cyclic and separating vector.
\item [(3)]  The right wandering
subspace of $K \cap L^2(M)$ (resp.\ $[K]_2$) in $L^2(M)$
has a cyclic and separating vector, when
$1 \leq p \leq 2$ (resp.\ $2 \leq p \leq \infty$).
\end{itemize} 
\end{proposition}

\begin{proof}  If $K = u [A]_p$ and $p < 2$, then
$K \cap L^2(M) = u H^2$ by an argument in the proof of (3)
of the Theorem.
  It follows that the right wandering subspace of $K \cap L^2(M)$ is
  $u [{\mathcal D}]_2$.  From this it is clear that the right wandering
subspace of $K$ is $u [{\mathcal D}]_p$.
These both have separating cyclic
vectors.  If $p > 2$ then this argument is easier (one
also uses the simple fact that $L^p(M) \cap L^2({\mathcal D})
= L^p({\mathcal D})$). 

To prove that (3) implies (1), note that by Corollary
\ref{stand} in the case $p = 2$, we have that
 $K \cap L^2(M)$ (resp.\ $[K]_2$) equals $u [A]_2$ for a unitary $u \in M$.
Thus $K = u [A]_p$ (if $p < 2$ use Saito's theorem 
mentioned in the introduction,
whereas if $p > 2$ use Lemma \ref{ll1}).
 
If (2) holds,
then by an adaption of an argument 
from \cite[p.\ 13]{JS} there exists 
an isometric ${\mathcal D}$-module isomorphism
  $\psi : L^p({\mathcal D}) \to W$.  (This is
a variant of the well known fact that 
a $W^*$-module over ${\mathcal D}$ with a cyclic separating vector
is unitarily isomorphic to ${\mathcal D}$; indeed the latter is
the case  $p = \infty$ of the assertion under discussion.)
Set $\psi(1) = u$
and set $v = u^* u \in L^{p/2}({\mathcal D})$.  
If $p < \infty$ we have $$\tau((d^* d)^{p/2}) = \Vert \psi(d) \Vert_p^p
=  \Vert u d \Vert_p^p =
\tau((d^* v d)^{p/2}) , \qquad d \in {\mathcal D} .$$
By Lemma \ref{istr}  with $p$ replaced by $p/2$,
 the last identity forces
$v = 1$, so that $u$ is unitary
(since we
are in a finite von Neumann algebra), and $W = u [{\mathcal D}]_p$.
Since $K_2^* W = (0)$, we have $K_2 = (0)$, and so $K = u H^p$.
A similar argument works if $p = \infty$; here 
in place of Lemma \ref{istr} one may use the fact that an 
isometric ${\mathcal D}$-module isomorphism between 
$C^*$-modules is unitary (see e.g.\ 8.1.5 in \cite{BLM}).  
  \end{proof}

The following is the analogue of Beurling's characterization of
weak* closed ideals of $H^\infty(\Ddb)$.  The last part uses
also the main theorem from \cite{N}, which is valid in $L^2$,
but which transfers easily using results above to $L^\infty$.    

\begin{corollary} \label{rid}  If  $A$ is maximal subdiagonal,
then the type 1 weak* closed right ideals of $A$
are precisely those right ideals of the form
$\oplus_i^{col} \, u_i A$, for partial isometries
$u_i \in A$ with mutually orthogonal ranges and $|u_i| \in {\mathcal D}.$
If the center of ${\mathcal D}$ is contained in the
center of ${\mathcal M}$, then one needs only one
partial isometry here.
\end{corollary}
 
{\bf Remark.}    As is proved in \cite{N} in the case $p = 2$,
a closed $A$-invariant subspace $K$ of $L^p(M)$ has type 1  
if and only if $[W M]_p = [K M]_p$ where $W$ is the 
right wandering subspace.  The one direction of this is obvious.
For the other note that if $K$ is not type 1 then by our main 
theorem there 
exists $\eta \in K$ with $\eta^* W = 0$.  Thus $\eta^* [W M]_p  = 0$
by e.g.\ (4.8) in \cite{JS}, which shows 
that $[W M]_p \neq [K M]_p$.

\begin{corollary} \label{pBN}
Suppose that $A$ is maximal subdiagonal,
and $f \in L^p(M)$, for $1 \leq p < \infty$.  If
$[f A]_p$ is a type 1 invariant subspace,
then $f = \sum_i u_i h_i$ (a norm convergent sum), where $u_i$ are  partial isometries in $[f A]_p
 \cap M$ with
$u_i^* u_i \in {\mathcal D}$, and $u_j^* u_i = 0$ if
$i \neq j$, and $h_i  \in [A]_p$ with $(u_i^* u_i) h_i = h_i$
and $u_i^* u_i \in [h_i A]_p$.
\end{corollary}

\begin{proof}   By the results above, $[fA]_p =
\oplus_i^{col} \, u_i [A]_p$, where
 $u_i$ are partial isometries in $M \cap [f A]_p$ of the correct form.
Moreover,  left multiplication
by the $u_i u_i^*$ are contractive projections from $[fA]_p$ onto
the summands $u_i [A]_p$.
Thus $f = \sum_i \, u_i h_i$ for $h_i \in (u_i^* u_i) [A]_p$.
We have $$u_i^* u_i \in (u_i^* u_i) [A]_p = u_i^* [f A]_p = (u_i^* u_i) [h_i A]_p
= [h_i A]_p .$$
This completes the proof.
\end{proof}

A similar result with almost identical proof holds if $p = \infty$,
interpreting closures and convergence in the weak* topology. 

The last result is a generalized `inner-outer' factorization.
The sum of products can be replaced by a single product if
the wandering subspace has a cyclic vector.   For example,
if $f \in L^p(M)$
and if $K = [f A]_p$ has a wandering subspace which has a 
cyclic and separating vector, then by Proposition \ref{triple}
we have $K = u H^p$ for a unitary $u \in M$.
Thus $f = u h$ for $h \in H^p$, and as in the last lines of the proof  of
Lemma  \ref{res}, this implies that
$[h A]_p = [A]_p$.  We take the 
latter condition as the definition of $h$ being {\em outer},
as in the classical case.

\bigskip

{\bf Closing remark:} 
In Sections 2 and 3, and in \cite{BL2},
we have been able to generalize almost the entire circle of equivalent
characterizations from \cite{SW}
of (at least weak* closed) weak*-Dirichlet algebras.
There are two items remaining in that list.
The first is known as the Gleason-Whitney theorem,
and we hope in the future to be able to demonstrate the equivalence of this
condition with the others, at least under some not too
unpalatable restriction on ${\mathcal D}$.
The second item is the condition that $A_\infty$ is
maximal subdiagonal.  Unfortunately we were unable to follow
the proof given in \cite{SW} for the latter equivalence,
nor have we been able to find this equivalence mentioned elsewhere in
the literature (without additional hypotheses).
One sufficient condition under which
$A_\infty$ being maximal subdiagonal implies that $A$ is 
maximal subdiagonal, is that the extension of $\Phi$ to $L^1(M)$ be continuous with
respect to the topology of \emph{convergence in measure}
(see
\cite{Terp,FK,Tak2,Nel} for details).  
However the latter does not hold for many interesting 
subdiagonal algebras.

\medskip

{\bf Acknowledgements.}  We thank Mike Marsalli, David Sherman,
and Dinesh Singh for valuable discussions.

\end{document}